\documentclass{amsart}

\usepackage[francais,english]{babel}
\usepackage[latin1]{inputenc}
\usepackage[T1]{fontenc}
\usepackage{ae,aecompl}

\title[Some relational structures and their associated algebras]
{Some relational structures with polynomial growth and their
  associated algebras}

\author{Maurice Pouzet and Nicolas M. Thiéry}
\address{UFR de Mathématiques,
  Université Claude-Bernard,
  $43$, Bd. du $11$ Novembre $1918$,
  $69622$ Villeurbanne, France,
  Fax 33 4 37 28 74 80}
\email{pouzet@univ-lyon1.fr, nthiery@users.sf.net}

\date{10th May 2005}

\newcommand{\TODO} [1]{}

\TODO{Florent: très dense; à délayer dans l'article final (plus tard)}
\TODO{Cite Pouzet-Sobrani? (non)}
\TODO{Finish the proof of the rationality of the Hilbert series}

\newcommand{\FIXME}[1]{}

\usepackage{amsfonts,amsthm,amssymb}
\usepackage{hyperref}

\usepackage{enumerate}
\usepackage{enumerate}

\newcommand{\isection}{{\downarrow}}
\newcommand{\sg}{{\mathfrak S}}
\newcommand{\lsg}{{\isection \sg}}
\newcommand{\bijection}{\hookrightarrow\!\!\!\!\!\rightarrow}

\newcommand{\dom}{{\operatorname{dom}}}
\newcommand{\im}{{\operatorname{im}}}
\newcommand{\rank}{{\operatorname{rank}}}
\newcommand{\id}{{\operatorname{id}}}
\newcommand{\support}{{\operatorname{support}}}
\newcommand{\gr}{{\operatorname{gr}}}

\newcommand{\aut}{{\operatorname{Aut}}}
\newcommand{\qsym}{{\operatorname{QSym}}}
\newcommand{\sym}{{\operatorname{Sym}}}
\newcommand{\lm}{{\operatorname{lm}}}
\newcommand{\hilbert}{{\mathcal H}}
\newcommand{\age}{{\mathcal A}}
\newcommand{\profile}{\varphi}
\newcommand{\lex}{{\operatorname{lex}}}

\newcommand{\R}{\mathbb{R}}
\newcommand{\K}{\mathbb{K}}
\newcommand{\N}{\mathbb{N}}
\newcommand{\C}{\mathbb{C}}
\newcommand{\Z}{\mathbb{Z}}
\newcommand{\Q}{\mathbb{Q}}

\newtheorem{theorem}{Theorem}[section]

\newtheorem{lemma}[theorem]{Lemma}

\newtheorem{proposition}[theorem]{Proposition} 
\newtheorem{corollary}[theorem]{Corollary} 
 
\theoremstyle{definition}

\newtheorem{problem}[theorem]{Problem}

\theoremstyle{remark}
\newtheorem{remark}[theorem]{Remark}
\newtheorem{remarks}[theorem]{Remarks}

\newtheorem{questions}[theorem]{Questions}
\newtheorem{example}[theorem]{Example}
\newtheorem{examples}[theorem]{Examples}

\def\centerpicture #1 by #2 (#3){\leavevmode
        \vbox to #2{
        \hrule width #1 height 0pt depth 0pt
        \vfill
        \special{pictfile #3}}}

\begin{document}

\maketitle

\hfill This paper is dedicated to Adriano Garsia at the occasion of
his 75th birthday

\begin{abstract}
  The \emph{profile} of a relational structure $R$ is the function
  $\profile_R$ which counts for every integer $n$ the number, possibly
  infinite, $\profile_R(n)$ of substructures of $R$ induced on the
  $n$-element subsets, isomorphic substructures being identified.
  Several graded algebras can be associated with $R$ in such a way
  that the profile of $R$ is simply the Hilbert function.  An example
  of such graded algebra is the \emph{age algebra} $\K.\age(R)$,
  introduced by P.~J.~Cameron. In this paper, we give a closer look at
  this association, particularly when the relational structure $R$
  decomposes into finitely many monomorphic components.  In this case,
  several well-studied graded commutative algebras (e.g.  the
  invariant ring of a finite permutation group, the ring of
  quasi-symmetric polynomials) are isomorphic to some $\K.\age(R)$.
  Also, $\profile_R$ is a quasi-polynomial, this supporting the
  conjecture that, with mild assumptions on $R$, $\profile_R$ is a
  quasi-polynomial when it is bounded by some polynomial.

\end{abstract}

\section{Presentation}

A \emph{relational structure} is a realization of a language whose
non-logical symbols are predicates, that is a pair $R:= (E,
(\rho_i)_{i \in I})$ made of a set $E$ and of a family of $n_i$-ary
relations $\rho_i$ on $E$. The family $\mu:=(n_i)_{i\in I}$ is the
\emph{signature} of $R$. The \emph{profile} of $R$ is the function
$\profile_R$ which counts for every integer $n$ the number
$\profile_R(n)$ of substructures of $R$ induced on the $n$-element
subsets, isomorphic substructures being identified.  Clearly, this
function only depends upon the set $\age(R)$ of finite substructures
of $R$ considered up to an isomorphism, a set introduced by R.~Fraïssé
under the name of \emph{age} of $R$ (see~\cite{Fraisse.TR.2000}). If
$I$ is finite $\profile_R(n)$ is necessarily finite. As we will see,
in order to capture examples coming from algebra and group theory, we
cannot preclude $I$ to be infinite. Since the profile is finite in
these examples, we will always make the assumption that
$\profile_R(n)$ is finite, no matter how large $I$ is.

A basic result about the behavior of the profile is this:
\begin{theorem}\label{increaseinfinite}
  If $R$ is a relational structure on an infinite set, then
  $\profile_R$ is non-decreasing.
\end{theorem}

This result was obtained in $1971$ by the first author (see Exercise~8
p.~113~\cite{Fraisse.CLM1.1971}). A proof based on linear algebra is
given in~\cite{Pouzet.1976}.

Provided that the relational structures satisfies some mild
conditions, there are jumps in the behavior of the profile:
\begin{theorem}\cite{Pouzet.TR.1978}\label{profilpouzet1}
  Let $R := (E, (\rho_i)_{i \in I})$ be a relational structure. The
  growth of $\profile_R$ is either polynomial or as fast as every
  polynomial provided that either the signature $\mu : = (n_i)_{i
    \in I}$ is bounded or the kernel $K(R)$ of $R$ is finite.
\end{theorem}
Note that a  map $\varphi : \N\rightarrow \N$ has {\it polynomial growth}, of {\it degree} $k$,   if $an^k\leq \varphi(n)\leq bn^k $ for some $a,b>0$  and $n$ large enough. The \emph{kernel} of $R$ is the set $K(R)$ of $x \in E$ such that
$\age(R_{\vert E \setminus \{x\}}) \neq \age(R)$.  Relations with
empty kernel are the \emph{inexhaustible} relations of R.~Fraïssé
(see~\cite{Fraisse.TR.2000}). We call \emph{almost inexhaustible}
those with finite kernel. The hypothesis about the kernel is not ad
hoc. As it turns out, if the growth of the profile of a relational
structure with a bounded signature is bounded by a polynomial then its
kernel is finite.  Some hypotheses on $R$ are needed, indeed \emph{for
  every increasing and unbounded map $\profile : \N \rightarrow \N$,
  there is a relational structure $R$ such that $\profile_R$ is
  unbounded and eventually bounded above by
  $\profile$}(cf.~\cite{Pouzet.RPE.1981}).

The consideration of examples suggested that in order to study the
profile of a relational structure $R$, the right object to consider is
rather the generating series
$$\hilbert_{\profile_R}:= \sum_{n=0}^\infty \profile_R(n)Z^n$$
This innocuous change leads immediately to several:
\begin{questions}\label{question}%
  \begin{enumerate}
  \item For which relational structures is the series
    $\hilbert_{\profile_R}$ a rational fraction? A rational fraction
    of the form
    \begin{displaymath}
      \frac{P(Z)}{(1-Z)(1-Z^{n_2})\cdots(1-Z^{n_k})}\ ,
    \end{displaymath}
    with $k\geq 1$, $1=n_1\leq n_2\leq \cdots \leq n_k$, $P(0)=1$, and
    $P\in \Z[Z]$?
  \item Is this the case of relational structures with bounded
    signature or finite kernel for which the profile is bounded by
    some polynomial?  When can $P$ be taken with nonnegative
    coefficients?
  \item For which relational structures is this series convergent?  Is
    this the case of relational structures $R$ whose age $\age(R)$ is
    well-quasi-ordered by embeddability?
  \end{enumerate}
\end{questions}
\begin{remark}\label{remark.quasipolynomial}
  When $\hilbert_{\profile_R}$ is a rational fraction of the form
  above then, for $n$ large enough, $\profile_R(n) $ is a
  \emph{quasi-polynomial} of degree $k'$, with $k'\leq k-1$, that is a
  polynomial $a_{k'}(n)n^{k'}+\cdots+ a_0(n)$ whose coefficients
  $a_{k'}(n), \dots, a_0(n)$ are periodic functions.  Since the
  profile is non-decreasing, it follows that $a_{k'}(n)$ is eventually
  constant. Hence the profile has polynomial growth:
  $\profile_R(n)\sim a n^{k'}$ for some non-negative real $a$.
\end{remark}
With the contribution of P.~J.~Cameron, this also gives links to some
quite venerable fields of mathematics. Indeed, P.~J.~Cameron~\cite
{Cameron.1997} associates to the age of $R$, $\age(R)$, its \emph{age
  algebra}, a graded commutative algebra $\K.\age(R)$ over a field
$\K$ of characteristic zero, and shows that the dimension of the
homogeneous component of degree $n$ of $\K.\age(R)$ is
$\profile_R(n)$, hence the generating series above is simply the
Hilbert series of $\K.\age(R)$.  Other graded algebras than the age
algebra enjoy this property. In any case, the association between
relational structures and graded algebras via the profile seems to be
an interesting topic.

The purpose of this paper is to document this association.  We do that
for a special case of relational structures that we introduce here for
the first time: \emph{those admitting a finite monomorphic
  decomposition}.  Despite the apparent simplicity of these relational
structures, the corresponding age algebras include familiar objects
like invariant rings of finite permutation groups. The profile of
these relational structures is a quasi-polynomial (cf. Theorem
\ref{theorem.quasipolynomial}). This supports the conjecture that the
profile of a relational structures with bounded signature or finite
kernel is a quasi-polynomial whenever the profile is bounded by some
polynomial (for more on the profile, see~\cite{Pouzet.PR.2002}).

We particularly study the special case of relational structures
associated with a permutation groupoid $G$ on a finite set $X$; as it
turns out, their age algebra is a subring of $\K[X]$, the
\emph{invariant ring} associated to $G$ (cf.
Theorem~\ref{polyalgebra}). This setting provides a close
generalization of invariant rings of permutation groups which includes
other famous algebras like quasi-symmetric polynomials and their
generalizations. We analyze in details which properties of invariant
rings of permutation groups carry over --- or not --- to permutation
groupoids (cf.  Propositions~\ref{proposition.groupoid.derivation}
and~\ref{proposition.groupoids.reynolds}, and
Theorems~\ref{theorem.finiteGeneration} and~\ref{theorem.SAGBI}). To
this end, we use in particular techniques
from~\cite{Garsia_Stanton.1984}.

A strong impulse to this research came from the paper by Garsia and
Wallach~\cite{Garsia_Wallach.2003}, which proves that, like invariant
rings of permutation groups, the rings of quasi-symmetric polynomials
are Cohen-Macaulay.  Indeed, a central long term goal is to
characterize those permutation groupoids whose invariant ring is
Cohen-Macaulay (cf. Problem~\ref{problem.CM}), this algebraic property
implying that the profile can be written as a quasi-polynomial whose
numerator has non-negative coefficients. As a first step in this
direction, we analyze several examples, showing in particular that not
all invariant rings of permutation groupoids are Cohen-Macaulay.

\section {Relational structures admitting a finite monomorphic decomposition}

A \emph{monomorphic decomposition} of a relational structure $R$ is a
partition $\mathcal P$ of $E$ into blocks such that for every integer
$n$, the induced structures on two $n$-elements subsets $A$ and $A'$
of $E$ are isomorphic whenever the intersections $ A\cap B$ and $
A'\cap B$ over each block $B$ of $\mathcal P$ have the same size.

\subsection{Some examples of relational structures admitting a finite
  monomorphic decomposition}.

\begin{example}
  Let $R:=(\Q, \leq, u_1,\dots,u_{k+1})$ where $\Q$ is the chain of
  rational numbers, $u_1,\dots,u_{k+1}$ are $k+1$ unary relations
  which divide $Q$ into $k+1$ intervals. Then $\profile_R(n) = {n+k
    \choose k}$ and $\hilbert_{\profile_R}=\frac{1}{(1-Z)^{k+1}}$.
\end{example}
\begin{example}
  A graph $G:=(V, \mathcal E)$ being considered as a binary
  irreflexive and symmetric relation, its profile $\profile_G$ is the
  function which counts, for each integer $n$, the number
  $\profile_G(n)$ of induced subgraphs on $n$ elements subsets of
  $V(G)$, isomorphic subgraphs counting for one.
  
  Trivially $\profile_G$ is constant, equal to $1$, if and only if $G$
  is a clique or an infinite independent set. A bit less trivial is
  the fact that $\profile_G$ is bounded if and only if $G$ is
  \emph{almost constant} in the sense of
  R.~Fraïssé~\cite{Fraisse.TR.2000}, that is there is a finite subset
  $F_G$ of vertices such that two pairs of vertices having the same
  intersection on $F_G$ are both edges or both non-edges.
\end{example}

\begin{example}
  Let $G$ be the direct sum $K_\omega\oplus K_\omega$ of two infinite
  cliques; then $\profile_G(n)=\lfloor \frac{n}{2} \rfloor +1$.
\end{example}

\begin{example}
  Let $G$ be the direct sum of $k+1$ many infinite cliques; then
  $\profile_G(n)=p_{k+1}(n) \simeq \frac{n^k} {(k+1)!k!}$ .
\end{example}

\begin{examples}
  Let $G$ be the direct sum $K_{(1, \omega)}\oplus \overline K_\omega$
  of an infinite wheel and an infinite independent set, or the direct
  sum $K_\omega\oplus \overline K_\omega$ of an infinite clique and an
  infinite set; then $\profile_G(n)=n$. Hence
  $\hilbert_{\profile_G}=1+\frac {Z}{(1-Z)^2}$, that we may write
  $\frac{1-Z-Z^2}{(1-Z)^2}$, as well as
  $\frac{1+Z^3}{(1-Z)(1-Z^2)}$.
\end{examples}
The first example above has three monomorphic components, one being
finite, whereas the second one has two component, both infinite;
still, the generating series coincide, with one representation as a
rational fraction with a numerator with some negative coefficient, and
another with all coefficients non-negative.

A more involved example is the following:
\begin{example}
  \TODO{Make a drawing, or check that this is just the same example as
    in \ref{example.nonCM.groupoid}, but with 3 variables. If yes,
    merge the examples.}%
  Let $R:=(E, (\rho, U_2,U_3) )$, where $E:=\N \times
  \{0,1,2,3\}$, $\rho:=\{((n,i),(m,j)): i=0, j\in \{1,2\}$ or
  $i=1,j=3\}$; $U_i:=\N\times \{i\}$ for $i\in \{0,1,2,3\}$. Then
  $R$ has four monomorphic components, namely $\N\times \{0\},\N
  \times \{1\},\N \times \{2\}, \N \times \{3\}$.  Let $S$ be the
  induced structure on four elements of the form $(x_i,i), i\in
  \{0,1,2,3\}$. A crucial property is that $S$ has only two
  non-trivial local isomorphisms, namely the map sending $(x_0,0)$
  onto $(x_1,1)$ and its inverse. From this follows that the induced
  substructures on two $n$-element subsets $E$ are isomorphic if
  either they have the same number of elements on each $\N\times
  \{i\}$ or one subset is included into $\N \times \{0\}$, the other
  into $\N\times \{1\}$. Hence, the generating series
  $\hilbert_{\profile_R}$ is $\frac{1}{(1-Z)^4}-\frac{Z}{1-Z}= \frac
  {1-Z+3Z^2-3Z^3+Z^4}{(1-Z)^4}$. We may write it
  $\hilbert_{\profile_R}=\frac {Q_1}{(1-Z)(1-Z^4)(1-Z^5)(1-Z^5)}$
  where
  $Q_1:=1+2Z+6Z^2+10Z^3+14Z^4+17Z^5+18Z^6+14Z^7+10Z^8+6Z^9+Z^{10}$,
  as well as $\hilbert_{\profile_R}=\frac {Q}{(1-Z)(1-Z^5)^3}$ where
  $Q_2:=1+2Z+6Z^2+10Z^3+15Z^4+18Z^5+22Z^6+18Z^7+15Z^8+10Z^9+6Z^{10}+
  Z^{12}+Z^{16}$.
\end{example}

Here is an example of relational structure $R$ such that
$\hilbert_{\profile_R}$ is a rational fraction, for which there is no
way of choosing the numerator with non-negative coefficients.

\begin{example}
  Let $R:=(E, (\rho_0, \rho_1, \rho_2) )$ be defined as follows.
  First, $E:= \N \times \{0,1,2,3\}\setminus \N^*\times \{0, 1\}$,
  that is $E$ is the union of two one-element sets, $E_0:=\{(0,0)\},
  E_1:=\{(0, 1)\}$, and of two infinite sets, $E_2:=\N \times \{2\},
  E_3:=\N \times \{3\}$.  Next $\rho_i:=E_i\times (E_2\cup E_3)$ for
  $i=0,1$ and $\rho_2:=E_0\times E_1\times (E_1 \cup E_2)$.  Then $R $
  has four monomorphic components, namely $E_0, E_1, E_2 E_3$. The
  crucial property is that the induced structures on two $n$-element
  subsets $A$, $B$ of $E$ are isomorphic if either $A$ and $B$
  includes $E_0\cup E_1$ and their respective traces on $ E_2$, $E_3$
  have the same size, or $A$ and $B$ include exactly one of the sets
  $E_0$ $E_1$ and not the other, or exclude both.  From this
  $\profile_R(0)=1$ and $\profile_R(n)=n+2$ for $n\geq 1$, hence
  $\hilbert_{\profile_R}= \frac {1+Z-Z^2}{(1-Z)^2}$. If
  $\hilbert_{\profile_R}= \frac {P}{(1-Z)(1-Z^k)}$ with $k\geq 2$,
  then $$P=1+2Z+ \sum_{j=2}^{k-2} Z^{j}-Z^{k+1}$$
  hence has a negative
  coefficient .
\end{example}

Another example with the same property, with two monomorphic component
which are both infinite.
\begin{example}\label{negative}
  Let $R:=(E, \mathcal H )$, where $E:=\N \times \{0,1\}$, $\mathcal
  H:=[\N \times \{0\}]^3\cup [\N \times \{1\}]^3$. Then $R$ has two
  monomorphic components, namely $\N \times \{0\}$ and $\N \times
  \{1\}$. Each type of $n$-element restriction has a representative
  made of a $m+k$ element subset of $\N \times \{0\}$ and of a
  $m$-element subset of $\N \times \{1\}$ such that $n=2m+ k$; these
  representatives are non-isomorphic, except if $n=2$ (in the later
  case, all $2$ -element restrictions are isomorphic, hence we may
  eliminate the representative corresponding to $m=1, k=0$). With this
  observation, a straightforward computation shows that $\varphi_R(0)=
  \varphi_R(1)=\varphi_R(2)=1$ and $\varphi_R(n)=\lfloor \frac{n}{2}
  \rfloor +1$ for $n\geq 3$. Hence the generating series
  $H_{\varphi_R}=\frac {1}{(1-x)(1-x^2)}- x^2=\frac
  {1-x^2+x^3+x^4-x^5}{(1-x)(1-x^2)}$.
  
  But, then $H_{\varphi_R}$ cannot be written as a quotient of the
  form $\frac{P}{(1-x)(1-x^k)}$ where $P$ is a polynomial with
  non-negative integer coefficients. Indeed, suppose, by
  contradiction, that $H_{\varphi_R}$ is of this form. We may suppose
  $k$ even (otherwise, multiply $P$ and $(1-x)(1-x^k)$ by $(1+x^k$).
  Set $k':=\frac{k}{2}$. Multiplying $1-x^2+x^3+x^4-x^5$ and
  $(1-x)(1-x^2)$ by $1+ x^2 + \cdots+ x^{2(k'-1)}$, we get
  $P=(1-x^2+x^3+x^4-x^5)(1+ x^2 + \cdots+ x^{2(k'-1)})$. Hence, the
  term of largest degree has a negative coefficient, a contradiction.

\end{example}

See also Example~\ref{example.nonCM.groupoid} for a last example with
three infinite monomorphic components and additional algebraic
structure which still has this property.

\subsubsection{Examples coming from group actions}

The \emph{orbital profile} of a permutation group $G$ acting on a set
$E$ is the function $\theta_G$ which counts for each integer $n$ the
number, possibly infinite, of orbits of the $n$-element subsets of
$E$.

As it is easy to see, orbital profiles are special cases of profiles.
Indeed, for every $G$ there is a relational structure such that $\aut
R = \overline{G}$ (the topological closure of $G$ in the symmetric
group $\mathfrak{G}(E)$, equipped with the topology induced by the
product topology on $E^E$, $E$ being equipped with the discrete
topology).  Groups for which the orbital profile takes only finite
values are said \emph{oligomorphic}, cf. P.~J.~Cameron
book~\cite{Cameron.1990}. These groups are quite common.  Indeed,
\emph{if $G$ is a group acting on a denumerable set $E$ and $R$ is a
  relational structure such that $\aut R = \overline{G}$ Then $G$ is
  oligomorphic if and only if the complete theory of $M$ is
  $\aleph_0$-categorical (Ryll-Nardzewski, 1959)}.

Even in the special case of groups, questions we ask in section 1 have
not been solved yet. For an example, let $G$ be a group acting on a
denumerable set $E$; if the orbital profile of $G$ is bounded above by
some polynomial, is the generating series of this profile a rational
fraction?  a rational fraction of the form given in Question 1? with a
numerator with non-negative coefficients?

Several examples come from relational structures which decompose into
finitely many monomorphic components.
\begin{example}
  Let $G$ be is the identity group on a $m$ element set $E$. Set $R:=
  (E, U_1,\dots,U_{m})$ where $U_1,\dots,U_{m}$ are $m$-unary
  relations defining the $m$ elements of $E$; then $\theta_G(n)=
  \profile_R (n)= {m \choose n}$.
\end{example}
\begin{example}
  Let $G := \aut\Q$, where $\Q$ is the chain of rational numbers. Then
  $\theta_G(n) = \profile_{\Q}(n)= 1$ for all $n$.
\end{example}

\begin{example}
  \label{example.permgroup}
  Let $G'$ be the wreath product $G':=G\wr \sg_\N$ of a permutation
  group $G$ acting on $\{1,\ldots, k\}$ and of $\sg_\N$, the symmetric
  group on $\N$. Looking at $G'$ as a permutation group acting on
  $E':=\{1,\ldots, k\}\times \N $, then $G'=\aut R'$ for some
  relational structure $R'$ on $E'$; moreover, for all $n$,
  $\theta_{G'}(n)=\profile_{R'}(n)$. Among the possible $R'$ take
  $R\wr \N:=(E', \equiv, (\overline \rho_i)_{i\in I})$ where $\equiv$
  is $\{((i, n),(j,m))\in E'^{2}: i=j\} $, $\overline \rho_i:=\{
  ((x_1, m_1),\dots,(x_{n_i}, m_{n_i})): (x_1,\dots, x_{n_i})\in
  \rho_i, ( m_1, \dots, m_{n_i})\in \N^{n_i} \}$, and $R:=(\{1,\dots,
  k\}, (\rho_i)_{i\in I})$ is a relational structure having signature
  $\mu:= (n_i)_{ i\in I}$ such that $\aut R= G$.  The relational
  structure $R\wr \N$ decomposes into $k$ monomorphic components,
  namely the equivalence classes of $\equiv$.
  
  As it turns out, $\hilbert_{\profile_{R\wr \N}}$ is the Hilbert
  series $\sum_{n=0}^\infty \dim \K[X]^G_n Z^n$ of the invariant ring
  $\K[X]^G$ of $G$ (that is the subring of the polynomials in the
  indeterminates $X:=(x_1,\dots,x_k)$ which are invariant under the
  action of $G$) (Cameron~\cite{Cameron.1990}).  As it is well known,
  this Hilbert series is a rational fraction of the form indicated in
  Question \ref{question}, where the coefficients of $P(Z)$ are
  non-negative.
\end{example}

\begin{problem}
  Find an example of a permutation group $G'$ acting on a set $E$ with
  no finite orbit, such that the orbital profile of $G'$ has
  polynomial growth, but the generating series is not the Hilbert
  series of the invariant ring $\K[X]^G$ of a permutation group $G$
  acting on a finite set $X$.
\end{problem}

\subsubsection{Quasi-symmetric polynomials and the like}

\begin{example}
  \label{example.qsym}
  Let $X_k:=(x_1,\ldots, x_k)$ be $k$ indeterminates and $n_1,\dots,
  n_l$ be a sequence of positive integers, $l\leq k$. The polynomial
  $$\sum_{1\leq i_1<\dots<i_l\leq k}x_{i_1}^{n_1}\dots x_{i_l}^{n_l}$$
  is a \emph{quasi-symmetric monomial} of degree
  $n:={n_1}+\dots+{n_l}$. The vector space spanned by the
  quasi-symmetric monomials forms the space $\qsym(X_k)$ of
  \emph{quasi-symmetric polynomials} as introduced by I.~Gessel. As in
  the example above, the Hilbert series of $\qsym(X_k)$ is defined as
  \begin{displaymath}
    \hilbert_{\qsym(X_k)}:= \sum_{n=0}^\infty \dim \qsym(X_k)_n Z^n\ .
  \end{displaymath}
  As shown by F.~Bergeron and C.~Reutenauer
  (cf.~\cite{Garsia_Wallach.2003}), this is a rational fraction of the
  form $\frac {P_k}{(1-Z)*(1-Z^2)*\dots(1-Z^k)}$ where the
  coefficients $P_k$ are non negative.  Let $R$ be the poset product
  of a $k$-element chain by a denumerable antichain. More formally,
  $R:= (E, \rho )$ where $E:= \{1,\ldots, k\}\times \N$ and $\rho:=
  \{((i, n), (j, m))\in E$ such that $i\leq j\}$.  Each isomorphic
  type of an $n$-element restriction may be identified to a
  quasi-symmetric polynomial, hence the generating series associated
  to the profile of $R$ is the Hilbert series defined above.
\end{example}

\begin{example}
  A relational structure $R:= (E, (\rho_i)_{i \in I})$ is
  \emph{categorical for its age} if every $R'$ having the same age as
  $R$ is isomorphic to $R$. It was proved in~\cite
  {Hodkinson_Macpherson.1988} that for relational structure with
  finite signature ($I$ finite) this happens just in case $E$ is
  countable and can be divided into finitely many blocks such that
  every permutation of $E$ which preserves each block is an
  automorphism of $R$. 
\end{example}

\subsection{Results and problems about relational structures admitting a
  finite monomorphic decomposition}

The following result motivates the introduction of the notion under
review.
\begin{theorem}\label{finimorphic}
  The profile of a relational structure $R$ is bounded by some integer
  if and only if $R$ has a monomorphic decomposition into finitely
  many blocks, at most one being infinite.
\end{theorem}
Relational structures satisfying the second condition of the above
sentence are the so-called \emph{almost-monomorphic relational
  structures} of R.~Fraïssé. Theorem~\ref{finimorphic} above was
proved in~\cite{Fraisse_Pouzet.1971} for finite signature and in~\cite
{Pouzet.RPE.1981}) for arbitrary signature by means of Ramsey theorem
and compactness theorem of first order logic.  From
Theorem~\ref{increaseinfinite} and Theorem~\ref{finimorphic}, it
follows that a relational structure $R$ has a monomorphic
decomposition into finitely many blocks, at most one being infinite if
and only if
\begin{displaymath}
  \hilbert_{\profile_R}=\frac {1+ b_1Z+ \dots +b_lZ^l} {1-Z}\ ,
\end{displaymath}
where $ b_1,\dots, b_l$ are non negative integers.

It is trivial that, if an infinite relational structure $R$ has a
monomorphic decomposition into finitely many blocks, whereof $k$ are
infinite, then the profile is bounded by some polynomial, whose degree
itself is bounded by $k-1$.

\TODO{Introduce first the "optimal monomorphic decomposition".}

\begin{theorem}\label{theorem.quasipolynomial}
  Let $R$ be an infinite relational structure $R$ with a monomorphic
  decomposition into finitely many blocks $(E_i, i\in X)$, $k$ of
  which being infinite. Then, the generating series
  $\hilbert_{\profile_R}$ is a rational fraction of the form:
  \TODO{Nicolas: Double check this!}
  \begin{displaymath}
    \frac{P(Z)}{(1-Z)(1-Z^2)\cdots(1-Z^k)}\ .
  \end{displaymath}
  In particular, remark~\ref{remark.quasipolynomial} applies.
\end{theorem}

To each subset $A$ of size $d$ of $E$, we associate the monomial
\begin{displaymath}
  x^{d(A)} := \prod_{i\in X} x_i^{d_i(A)}\ ,
\end{displaymath}
where $d_i(A)=|A\cap E_i|$ for all $i$ in $X$. Obviously, $A$ is
isomorphic to $B$ whenever $x^{d(A)}=x^{d(B)}$. The shape of a
monomial $x^d=\prod x_i^{d_i}$ is the partition obtained by sorting
decreasingly $(d_i, i\in X)$. We define a total order on monomials by
comparing their shape w.r.t. the degree reverse lexicographic order,
and breaking ties by the usual lexicographic order on monomials w.r.t.
some arbitrary fixed order on $X$. To each orbit of sets, we associate
the unique maximal monomial $\lm(A)$, where $A$ ranges through the
orbit; we call this monomial \emph{leading monomial}. To prove the
theorem, we essentially endow the set of leading monomials with an
ideal structure in some appropriate polynomial ring. This is
reminiscent of the chain-product technique as defined in
Subsection~\ref{invringper}.  The key property of leading monomials is
this:

\begin{lemma}\label{lemma.addlayer}
  Let $m$ be a leading monomial, and $S\subset X$ a layer of $m$.
  Then, either $d_i=|E_i|$ for some $i$ in $S$, or $m x_S$ is again a
  leading monomial.
\end{lemma}

The proof of this result relies on
Proposition~\ref{homogeneouscomponent} below for which we introduce
the following definition. Let $R$ be a relational structure on $E$; a
subset $B$ of $E$ is a \emph{monomorphic part} of $R$ if for every
integer $n$ and every pair $A, A' $ of $n$-element subsets of $E$ the
induced structures on $A$ and $A'$ are isomorphic whenever $A\setminus
B=A'\setminus B$. The following lemma, given without proof, rassembles
the main properties of monomorphic parts.
\begin{lemma}\label{lemmahomogeneous}
  \begin{enumerate}[(i)]
  \item The emptyset and the one element subsets of $E$ are
    monomorphic parts of $R$;
  \item If $B$ is a monomorphic part of $R$ then every subset of $B$
    too;
  \item Let $B$ and $B' $ be two monomorphic parts of $R$; if $B$ and
    $B'$ intersect, then $B\cup B'$ is a monomorphic part of $R$;
  \item Let $\mathcal B$ be a family of monomorphic parts of $R$; if
    $\mathcal B$ is up-directed (that is the union of two members of
    $\mathcal B$ is contained into a third one), then their union $B:=
    \bigcup \mathcal B$ is a monomorphic part of $R$.
  \end{enumerate}
\end{lemma}
Let $x\in E$, let $R(x)$ be the set-union of all the monomorphic parts
of $R$ containing $x$.  By $(i)$ of Lemma~\ref{lemmahomogeneous} this
set contains $x$ and by $(iii)$ and $(iv)$ this is a monomorphic part,
thus the largest monomorphic part of $R$ containing $x$.

\begin {proposition}\label{homogeneouscomponent}
  The largest monomorphic parts form a monomorphic decomposition of
  $R$ off which every monomorphic decomposition of $R$ is a
  refinement.
\end{proposition}
\begin{proof}[Proof of Lemma \ref {lemma.addlayer}]
  Let $e:= \vert X\vert$, $\overline d:=(d_{i_1},\dots,d_{i_e})$ be
  the shape of $m$ sorted decreasingly and $s:= \vert S\vert$. Suppose
  that $d_i<|E_i|$ for every $i$ in $S$. Let $A$, $B$ , $B'$ be
  subsets of $E$ such that $x^{d(A)}= m$, $x^{d(B)}= mx_{S}$, and
  $m':= x^{d(B')}$ is the leading monomial in the orbit of $B$ and let
  $R_A$, $R_B$, $R_{B'}$ be the corresponding induced structures.
  
  Clearly, the shape of $mx_{S}$ is $\overline
  {d_1}:=(d_{i_1}+1,\dots, d_{i_s}+1, d_{i_{s+1}}, \dots, d_{i_e})$.
  Let $\overline {d'}:= (d'_{i'_1},\dots, d'_{i'_s}, d'_{i'_{s+1}},
  \dots, d'_{i'_e})$ be the shape of $B'$. Our first goal is to prove
  that these two shapes are the same.
  
  {\bf Claim 1} {\it $d_{i_p}=d'_{i'_p}$ for all $p> s$}.
  
  {\bf Proof of Claim 1} Suppose this does not hold. Let $p$ be the
  largest such that $d_{i_p}\not = d'_{i'_p}$.  Since, by definition,
  we have $\overline {d'}\geq \overline {d_1}$ it follows that
  $d'_{i'_p}<d_{i_p}$; thus, $\overline {d'}>\overline {d}$. However,
  since $R_{B'}$ contains a copy of $R_{A}$ we have $\overline
  {d'}\leq \overline {d}$, a contradiction.  \endproof
  
  Set $U:= \bigcup \{E_{i}\cap B : i\not \in S\}$, $S':= \{i'_1,
  \dots, i'_s\}$, $U':= \bigcup \{E_{i}\cap B': i\not \in S'\}$. Let
  $\varphi$ be an isomorphism from $R_B$ onto $R_{B'}$.
  
  {\bf Claim 2} {\it $U$ is the set of $x\in B$ such that the induced
    structure $R_{B\setminus \{x\}}$ on $B\setminus \{x\}$ contains no
    copy of $R_A$. Moreover $\varphi$ transforms $U$ into $U'$}.
 
  {\bf Proof of Claim 2} From the definition of $U$, $R_{B\setminus
    \{x\}}$ contains a copy of $R_A$ for every element $x \in
  B\setminus U$. Conversely, let $x\in U$ and let $\overline{d''}$ be
  the shape of $B\setminus \{x\}$. Clearly, for the largest $p$ such
  that $d''_{i''_p}\not = d_{i_p}$ we have $p>s$.  Hence $\overline
  d'' > \overline d$, thus $R_{B\setminus \{x\}}$ cannot contains a
  copy of $R_A$. This proves the first part of Claim 2.
  
  Since, from Claim 1, $d_{i_p}=d'_{i'_p}$ for all $p\geq s$, the same
  argument show that if $x'\in U'$ then $R_{B'\setminus \{x'\}}$
  cannot contains a copy of $R_A$. Since from Claim 1, $U$ and $U'$
  have the same size we get that $U'$ is the set of $x' \in B'$ such
  that $R_{B'\setminus \{x'\}}$, contains no copy of $R_A$. The second
  part of Claim 2 follows immediately.  \endproof
  
  {\bf Claim 3} {\it Let $i\not \in S$ and $j\in S$ then every
    monomorphic part containing $E_i\cap B$ is disjoint from $E_j\cap
    B$}.
  
  {\bf Proof of Claim 3} According to Claim 2 a monomorphic part
  containing $E_i\cap B$ must be disjoint from $U$. \endproof

  {\bf Claim 4} {For each $i\in S$, $E_i\cap B$ is a largest
    monomorphic part of $R_{B}$}.
  
  {\bf Proof of Claim 4} Suppose not. Then this largest monomorphic
  part, say $C$, contains some other $E_j\cap B$. From Claim 3, $j\in
  S$. It follows that all induced substructures on $C\setminus
  \{x,y\}$, where $\{x,y\}$ is a pair of distinct elements of $C$, are
  isomorphic.  Suppose $d_i\geq d_j$.  Since the shape of $A$ is
  maximal then for $x,y\in E_j$ the induced structure does not contain
  a copy of $R_A$. But if $x\in E_i$ and $y\in E_j$ then trivially the
  induced structure contains a copy of $R_A$.  A contradiction.
  \endproof
  
  {\bf Claim 5} {\it $\varphi$ transforms $(E_i\cap B, i\in S)$ into
    $(E_i\cap B', i\in S')$}
 
  {\bf Proof of Claim 5}.  The $E_i\cap B$'s for $ i\in S$ are the
  largest monomorphic parts of $R_{B\setminus A}$. Via $\varphi$ there
  are transformed into the $s$ largest monomorphic parts of
  $R_{B'\setminus U'}$.  Since $(E_i\cap B', i\in S')$ is a
  decomposition of $R_{B'\setminus U'}$ into $s$ monomorphic parts,
  this decomposition coincides with this decomposition into largest
  parts.\endproof
 
  From Claim 1 and Claim 5, we have $\overline {d'}= \overline {d_1}$.
  Suppose that $m'>mx_S$. Let $T$ be a transversal of the $E_{i}\cap
  B$'s for $ i\in S$.  Then, from Claim 5, $T':=\varphi (T)$ is a
  transversal of the $E_{i'}\cap B'$'s for $ i'\in S'$. Let $m_{T}$,
  resp. $m_{T'}$, be the monomial associated with $B\setminus T$,
  resp.  $B'\setminus T'$. We have $m'_{T'}>m_{T}$.  Since $m_{T}=m$
  and $B'\setminus T'$ is in the orbit of $B$, we get a contradiction.
\end{proof}

\begin{proof}[Proof of theorem~\ref{theorem.quasipolynomial}]
  Fix a chain $C=(\emptyset \subsetneq S_1\subsetneq\dots\subsetneq
  S_r \subset X)$ of non empty subsets of $X$.  Let $\lm_C$ be the set
  of leading monomials with chain support $C$.  The plan is
  essentially to realize $\lm_C$ as the linear basis of some ideal of
  a polynomial ring, so that the generating series of $\lm_C$ is
  realized as an Hilbert series.
  
  Consider the polynomial ring $\K[S_1,\dots,S_l]$, with its natural
  embedding in $\K[X]$ by $S_j\mapsto \prod_{i\in S_j} x_i$. Let $I$
  be the subspace spanned by the monomials $m=S_1^{r_1}\dots
  S_l^{r_l}$ such that $d_i(m)>|E_i|$ for some $i$; it is obviously a
  monomial ideal. When all monomorphic components are infinite, $I$ is
  the trivial ideal $\{0\}$. Consider the subspace $\K.\lm_C$ of
  $\K[S_1,\dots,S_l]$ spanned by the monomials in $\lm_C$.
  Lemma~\ref{lemma.addlayer} exactly states that $J=\K.\lm_C\oplus I$
  is in fact a monomial ideal of $\K[S_1,\dots,S_l]$. Both $I$ and $J$
  have finite free resolution as modules over $\K[S_1,\dots,S_l]$, so
  that their Hilbert series are rational fractions of the
  form\TODO{citation}:
  \begin{displaymath}
    \frac{P}{ (1-Z^{|S_1|}) \cdots (1-Z^{|S_l|}) }\ .
  \end{displaymath}
  Hence, the same hold for
  $\hilbert_{\K.\lm_C}=\hilbert_J-\hilbert_I$. Furthermore, whenever
  $S_j$ contains $i$ with $|E_i|<\infty$, the denominator
  $(1-Z^{|S_l|})$ can be canceled out in $\hilbert_{\K.\lm_C}$. The
  remaining denominator divides $(1-Z) \cdots (1-Z^k)$.
  \TODO{Improve the argument.}

  By summing up those Hilbert series $\hilbert_{\K.\lm_C}$ over all
  chains $C$ of subsets of $X$, we get the generating series of all
  the leading monomials, that is the profile of $R$. Hence, this
  profile is a rational fraction of the form:
  \begin{displaymath}
    \hilbert_{\K.\lm_C} =
    \frac{P}{ (1-Z) \cdots (1-Z^k) }\ .
  \end{displaymath}
\end{proof}

\begin{remark}
  As Examples~\ref{negative} and~\ref{example.nonCM.groupoid}
  illustrates, it is not true that if all blocks of a monomorphic
  decomposition of $R$ are infinite, then the numerator $P$ in the
  above fraction can be choosen with non-negative coefficients.
\end{remark}
We do not know for which relational structures having a finite monomorphic decomposition the numerator $P$ can be choosen with non-negative coefficients.  
A possible approach is to look for some
sensible Cohen-Macaulay graded algebra whose Hilbert series is
$\hilbert_{\profile_R}$ (by proposition~4 of~\cite{Bogvad_Meyer.2003}
such a Cohen-Macaulay algebra always exists as soon as $P$ has
non-negative coefficients). This is one of our motivations for the
upcoming study of the \emph{age algebras}.

\section{The age algebra of a relational structure}

\subsection{The set-algebra}

Let $E$ be a set and let $[E]^{<\omega}$ be the set of finite subsets
of $E$ (including the empty set). Let $\K$ be a field%
, and $\K^{[E]^{<\omega}}$
be the set of maps $f:[E]^{<\omega}\rightarrow \K$.  Endowed with the
usual addition and scalar multiplication of maps, $\K^{[E]^{<\omega}}$
is a $\K$-vector space. Let $f,g\in \K^{[E]^{<\omega}}$; according to
Cameron, we set:
\begin{displaymath}
  f g(P) := \sum_{M \in [P]^{<\omega}} f(P)g(P\setminus M)\ ,
\end{displaymath}
for all $P\in [E]^{<\omega}$~\cite{Cameron.1997}. With this operation
added, $\K^{[E]^{<\omega}}$ becomes a ring. This ring is commutative
and has a unit: denoted by $1$, this is the map taking the value $1$
on the empty set and the value $0$ everywhere else.\\
Let $\equiv$ be an equivalence relation on $[E]^{<\omega}$. A map
$f:[E]^{<\omega}\rightarrow \K$ is \emph{$\equiv$-invariant} or,
briefly, \emph{invariant} if $f$ is constant on each equivalence
class.  Invariant maps form a subspace of the vector space
$\K^{[E]^{<\omega}}$. We give a condition below which insures that
they form a subalgebra too.
\begin{lemma}\label{hereditary1}
  Let $\equiv$ be an equivalence relation on $[E]^{<\omega}$ and
  $D,D'\in [E]^{<\omega}$.  Then the following properties are
  equivalent:
  \begin{enumerate}[(i)]
  \item There exists some bijective map $f:D\bijection D'$ such that
    $D\setminus \{x\} \equiv D'\setminus \{f(x)\}$ for every $x\in D$;
  \item
    \begin{itemize}
    \item [1)]$\vert D\vert= \vert D'\vert=d$ for some $d$;
    \item [2)] $\vert \{X\in [D]^{d-1}:X\equiv B\}\vert = \vert \{X\in
      [D']^{d-1}: X\equiv B\}\vert$ for every $B\subseteq E$.
    \end{itemize}
  \end{enumerate} 
\end{lemma}
An equivalence relation on $[E]^{<\omega}$ is \emph{hereditary} if
every pair $D,D'$ of equivalent elements satisfies one of the two
equivalent conditions of Lemma~\ref{hereditary1}.
\begin{remark}
  Hereditary equivalences are introduced
  in~\cite{Pouzet_Rosenberg.1986} with Condition (ii)~2) of
  Lemma~\ref{hereditary1} replaced by the condition:
  \begin{displaymath}
    \vert \{X\subseteq D:X\equiv B\}\vert = \vert \{X\subseteq D':
    X\equiv B\}\vert \; \mbox{ for every}\; B\subseteq E\ .
  \end{displaymath}
  It follows from the next Lemma that this condition is not stronger.
\end{remark}

Let $\equiv$ be an equivalence relation on $[E]^{<\omega}$. We denote
by $[E]^{<\omega}_{/\equiv}$ the set of equivalence classes.
Let $a,b,c\in [E]^{<\omega}_{/\equiv}$ and $D\in [E]^{<\omega}$. Set
\begin{displaymath}
  \chi_{a,b,c}(D):=\vert\{(A,B)\in a\times b: A\cup B=C, C\subseteq D,
  C\in c\}\vert\ .
\end{displaymath}
If all subsets of $E$ belonging to some equivalence class $a$ have the
same size, we denote by $\vert a\vert $ this common size.

\begin{lemma}\label {hereditary2}
  If $\equiv$ is an hereditary equivalence relation on $[E]^{<\omega}$
  then
  \begin{equation}
    \chi_{a,b,c}(D):=\chi_{a,b,c}(D')\; \text{whenever }\; D\equiv D'.
  \end{equation}
\end{lemma}

\begin{proposition}\label{hereditary 3}
  Let $\equiv$ be an hereditary equivalence relation on
  $[E]^{<\omega}$.  Then the product of two invariant maps is
  invariant.
\end{proposition}

\subsection{The age algebra}

Let $R$ be a relational structure with domain $E$.  Set $F\equiv F'$
for $F,F'\in [E]^{<\omega}$ if the restrictions $R\restriction_{F}$
and $R\restriction_{F'}$ are isomorphic. The resulting equivalence on
$[E]^{<\omega}$ is hereditary, hence the set of invariant maps
$f:[E]^{<\omega}\rightarrow \K$ form a subalgebra of
$\K^{[E]^{<\omega}}$. Let $\K.\age(R)$ be the subset made of the
invariant maps which are everywhere zero except on a finite number of
equivalence classes. Then $\K.\age(R)$ forms an algebra, the \emph{age
  algebra} of Cameron.

\section{Invariant rings of permutation groupoids}

Let $\K$ be a field of characteristic $0$. In this section, we study
in more details a specific class of age algebras which can be realized
as graded subrings of polynomial rings $\K[X]$ that we call
\emph{invariant rings of permutation groupoids}.  This class extends
the class of invariant rings of permutation groups
(Example~\ref{example.permgroup}), and contains other interesting
examples like the rings of quasi-symmetric polynomials
(Example~\ref{example.qsym}).  Our long-term motivations are twofold.
On one hand, relate, in this simpler yet rich setting, the properties
of the profile to algebraic properties of the invariant ring. In
particular, find conditions under which the invariant ring is
Cohen-Macaulay.  On the other hand, generalize the theory, algorithms,
and techniques of invariant rings of permutation groups to a larger
class of subrings of $\K[X]$. In particular, find new properties of
the ring of quasi-symmetric polynomials, one specific goal being to
find a simpler proof that this ring is Cohen-Macaulay.


\subsection{Permutation groupoids}

Let $X$ be a finite set. A \emph{local bijection} of $X$ is a
bijective function $f: \dom f \bijection \im f$ whose domain $\dom f$
and image $\im f$ are subsets of $X$. The \emph{rank} of $f$ is the
size of its domain, so that $f$ is a permutation of $X$ if it is of
maximal rank $|X|$. The inverse $f^{-1}$ of a local bijection $f$, its
restriction $f_{\restriction X'}: X'\bijection f(X')$ to a subset $X'$
of $\dom f$, and the composition $f\circ g$ of two local bijections
$f$ and $g$ such that $\im g=\dom f$ are defined in the natural way.
A set $G$ of local bijections of $X$ is called a \emph{permutation
  groupoid} if it contains the identity and is stable by restriction,
inverse, and composition. It can be seen as a category: the objects
are the subsets of $X$ and the morphisms $G(A,B)$ from $A$ to $B$ are
the local bijections $f: A\bijection B$ in $G$. Those morphisms are by
definition isomorphisms, and $G$ satisfies the usual
\href{http://en.wikipedia.org/wiki/Groupoid}{\emph{groupoid} axioms}.

The \emph{underlying permutation group} is the subset $G(X,X)$ of all
permutations in $G$; those are exactly the invertible elements w.r.t.
the composition product.

\begin{examples}
  The set $\lsg(X)$ of all local bijections of $X$ is a permutation
  groupoid.
  
  The closure $\isection G$ of a permutation group $G$ by restriction
  is a permutation groupoid. In the following, we say that $\isection
  G$ \emph{comes from the permutation group} $G$.
  
  Let $X:=\{1,\dots,n\}$. The set $G$ of strictly increasing local
  bijections of $X$ forms a permutation groupoid.  Obviously $G$ does
  not come from a permutation group since its underlying permutation
  group is reduced to the identity.
  
  Let $R$ be a relational structure on $X$. The local isomorphisms of
  $R$ form a permutation groupoid. Its underlying permutation group is
  the automorphism group of $R$. Typically, the previous example is
  obtained by taking as relational structure $R$ the chain
  $1<2<\dots<n$. Also, $\lsg(X)$ is obtained by taking the trivial
  relational structure on $X$. In fact, \TODO{prove or reference it!}
  {\it any permutation groupoid $G$ can be obtained from a suitable
  relational structure $R_G$ on $X$} (recall that $X$ is finite!).
\end{examples}




\subsubsection{The invariant ring of a permutation groupoid}\label{invringper}

Let $G$ be a permutation groupoid acting on a finite set $X$, and
$\K[X]$ be the polynomial ring whose variables $x_i$ are indexed by the
elements $i$ of $X$. Given a local function $f$ of $G$, and a monomial
$x^d := \prod_{i\in X} x_i^{d_i}$ whose support $\support(x^d)$ is
contained in $\dom f$, we set
\begin{displaymath}
  f.x^d := \prod_{i\in X, d_i>0} x_{f(i)}^{d_i}\ ,
\end{displaymath}
generalizing the usual action of a permutation on a monomial.  This
partial action of $G$ on monomials does not extend to a global action
of $G$. Still, notions like \emph{$G$-isomorphic monomials} and
\emph{$G$-orbits} are well defined. The \emph{orbit sum} $o(x^d)$ of a
monomial $x^d$ is the sum of all the monomials in its orbit.

Our object of study is the \emph{invariant ring} $\K[X]^G$ of $G$,
which is defined as the linear subspace of $\K[X]$ spanned by the
orbitsums of all monomials. 

\begin{examples}
  Let $G$ be a permutation group. Then, $\K[X]^{\isection G}$ is the
  usual invariant ring of $G$.
  
  Let $G$ be the permutation groupoid of the strictly increasing local
  bijections of $\{1,\dots, n\}$. Then $\K[X]^G$ is the ring
  $\qsym(X_n)$ of quasi-symmetric polynomials on the ordered alphabet
  $X_n:=(x_1,\dots,x_n)$.
  
  Taking the same groupoid $G$ as in the previous example, and letting
  it act naturally on respectively pairs, couples, $k$-subsets, or
  $k$-tuples of elements of $\{1,\dots,n\}$, yields respectively the
  (un)oriented (hyper)graph quasi-symmetric polynomials
  of~\cite{Novelli_Thibon_Thiery.2004}.
\end{examples}

\begin{remarks}
  The orbitsums form a linear basis of $\K[X]^G$.
  
  It is not obvious from the definition that $\K[X]^G$ is indeed a
  graded algebra.  In the following we prove this by making $G$ into a
  monoid and the action of $G$ on polynomials into a multiplicative
  linear representation of $G$.
  An other way is to encode $G$ by some relational structure $R_G$ on
  $X$ and, as in Example~\ref{example.permgroup}, to define a
  relational structure $R_G \wr \N$ on $E:= X\times \N$ with
  monomorphic components $E_i:=\{i\}\times \N$ for $i\in X$.  Let
  $\phi: \K[X] \hookrightarrow \K^{[X\times \N]^{<\omega}}$ defined by
  setting $\phi(x^d ):= d! \chi_{O_{\mathfrak G}(x^d)}$, where $d!:=
  \prod_{i\in X} d_i!$, and $\chi_{O_{\mathfrak G}(x^d)}$ is the
  characteristic function of $O_{\mathfrak G}(x^d):=\{A\subset
  X\times\N,\, |A\cap E_i|=d_i,\forall i\in X\}$. Once $\K^{[X\times
    \N]^{<\omega}}$ is equipped with its set-algebra stucture, $\phi$
  is a morphism of algebras. Applying Proposition~\ref{hereditary 3},
  we get:
\end{remarks}

\begin{theorem}\label{polyalgebra}
  The invariant ring $\K[X]^G$ is isomorphic via $\phi$ to the age
  algebra $\K.\age(R_G \wr \N)$. In particular, the generating series
  of the orbits, the generating series of the profile of $R_G\wr\N$,
  and the Hilbert series of $\K[X]^G$ coincide.
\end{theorem}

\subsubsection{Restrictions of permutation groupoids}

The \emph{restriction} $G_{\restriction X'}$ of a permutation groupoid
$G$ to a subset $X'$ is the set of all local functions $f$ in $G$ such
that $\dom f\subset X'$ and $\im f \subset X'$, which is again a
permutation groupoid.  Furthermore, the orbits of
monomials in $\K[X']$ are unchanged by this restriction. In particular, the
invariant ring of $G_{\restriction X'}$ is simply the quotient of the
invariant ring of $G$ obtained by killing all the variables $x_i$ with
$i\not\in X'$.
This simple fact is one of the points of considering permutation
groupoids instead of just permutation groups (for which the
restriction to a subset is not clearly defined). This may indeed give
opportunities for induction techniques on the size of the underlying
set.
\begin{proposition}
  Any permutation groupoid comes from the restriction of a permutation
  group of some superset. However, this superset may need to be
  infinite.
\end{proposition}

\begin{examples}
  \label{examples.restrictions}
  (a) The permutation groupoid on $\{1,2,3\}$ generated by the
  rank $1$ local bijection $1\mapsto 2$ is the restriction of the
  permutation group on $\{1,2,3,4\}$ generated by the permutation
  $(1,2)(3,4)$.

  (b) The local automorphism permutation groupoid of the chain $a<b$
  is the restriction of the cyclic group $C_3$ on $\{a,b,c\}$.
  
  (c) Consider a relational structure $R$ such that there exists three
  elements $a,b,c$ and a binary relation $<$ which restricts on
  $\{a,b,c\}$ to the chain $a<b<c$. Typically, $R$ is a chain of
  length at least $3$ (giving $\qsym(X)$ as invariants) or a poset of
  height at least $3$. Then, there exists no relational structure
  $\overline R$ on a finite superset where all local isomorphisms
  extend to global isomorphisms.

  \FIXME{Uncomment the 4-lines proof?}
\end{examples}

\subsubsection{The monoid of a permutation groupoid}

The goal is now to turn $G$ into a monoid, and to make the partial
actions of $G$ into a linear representations of this monoid. The
composition of two local functions $f$ and $g$ can be extended when
$\im f\ne \im g$ by setting it to the local function with the largest
domain on which $f(g(x))$ is well-defined:
\begin{displaymath}
  f\circ g:
  \begin{cases}
    g^{-1}(\im g \cap \dom f) &\bijection f(\im g \cap \dom g)\\
    x \mapsto f(g(x))
  \end{cases}\ .
\end{displaymath}
With this composition product, $G$ turns into a monoid whose unit is
the identity of $X$. %
\TODO{properties: the local identities $\id_{A}$ on subsets of $A$ of
  $X$ are idempotent, $\id_{\{\}}$ is "absorbant", ...}%
Now we can extend the partially defined action of local bijections on
monomials into a linear action on polynomials by setting:
\begin{displaymath}
  f.x^d :=
  \begin{cases}
    \displaystyle
    \prod_{i\in X, d_i>0} x_{f(i)}^{d_i} & \text{ if $\support(x^d)
      \subset \dom f$,}\\
    0 & \text{ otherwise.}
  \end{cases}
\end{displaymath}

We leave it as exercise to check that this defines a linear
representation of the monoid $G$, which is \emph{multiplicative}: for
any $f$ in $G$ and $P$ and $Q$ in $\K[X]$,
\begin{displaymath}
  f.(PQ) = (f.P) (f.Q)\ .
\end{displaymath}

\begin{corollary}
  The invariant ring $\K[X]^G$ is, as its name suggests, indeed a
  ring.
\end{corollary}
\begin{proof}
  Consider a product of two orbitsums $o(m_1)o(m_2)$, and take two
  isomorphic monomials $m$ and $f.m$, $f\in G$. Whenever $m$ occurs as
  a product $m=m'_1m'_2$, $m'_1\in G.m_1$, $m'_2\in G.m_2$, the
  monomial $f.m$ occurs simultaneously as the product
  $f.m=f.(m'_1m'_2)=f.m'_1f.m'_2$, and reciprocally. Hence $m$ and
  $m'$ occur with the same coefficient in $o(m_1)o(m_2)$.
\end{proof}

Note that the monoid algebra of $G$ is isomorphic to its groupoid
algebra $\K.G$ which is semi-simple. This linear representation of $G$
extends into a linear representation of $\K.G$.


\subsubsection{Groupoid and monoid algebra of a permutation groupoid}

Let $G$ be a permutation groupoid, and $\K$ a field (of characteristic
zero; typically $\K=\R$ or $\K=\C$). By definition, its \emph{groupoid
  algebra} $\K.G$ is the $\K$-vector space whose basis $\{\gr f, f\in
G\}$ is indexed by the elements $f$ of $G$, and whose product is given
by:
\begin{displaymath}
  \gr f \ \gr g =
  \begin{cases}
    f\circ g & \text{ if $\im g=\dom f$, }\\
    0 & \text{otherwise.}
  \end{cases}
\end{displaymath}
We call $\{\gr f, f\in G\}$ the \emph{graded basis} of $\K.G$.
Similarly the \emph{monoid algebra} of $G$ is defined as the
$\K$-vector with basis $\{f, f\in G\}$ equipped with the extended
composition product $\circ$.

\begin{remarks}
  As the notation suggests, the groupoid and the monoid algebra of $G$
  are isomorphic.

  The groupoid algebra $\K.G$ is semi-simple, and decomposes as a
  direct sum of non-unitary algebras:
  \begin{displaymath}
    \K.G = \sum_{k=0}^n \K.G_k\ ,
  \end{displaymath}
  where $G_k := \{ f\in G, \rank f=k\}$.
\end{remarks}

\begin{proof}
  The isomorphism from the monoid algebra to the groupoid algebra is
  given by:
  \begin{displaymath}
    f \mapsto \sum_{A\subset \dom f} \gr f_{\restriction A}\ .
  \end{displaymath}
  The inverse isomorphism is obtained by Möbius inversion:
  \begin{displaymath}
    \gr f \mapsto \sum_{A\subset \dom f} (-1)^{|\dom f|-|A|}
    f_{\restriction A}\ .
  \end{displaymath}
  Checking the compatibility with the product rule is straightforward.

  The semi-simplicity is a general property of groupoid algebras,
  which is easily checked using Dickson's Lemma.
\end{proof}
  
  


The linear representations of the monoid $G$ on
polynomials extend directly to linear representations of its
algebra $\K.G$. In particular, this defines the actions of the graded
basis. Its characteristic is that $\gr f$ kills all monomials whose
support is not exactly $\dom f$, whereas $f$ kills only those monomials
whose support is not contained in $\dom f$.

Important note: the action of $\gr f$ is \emph{not} multiplicative on
polynomials! Take for example $f:=\id_{\{1,2\}}$, $P:=x_1$ and
$Q:=x_2$. This is in fact the main reason for considering the monoid
algebra and not just only the groupoid algebra.





\subsection{Invariants of permutation groupoids}

In this section, we review which properties of invariants of
permutation groups extend to permutation groupoids.

\subsubsection{The Reynolds operator}

The first essential feature of invariant rings is the so-called
\emph{Reynolds operator}, which is a projector on the invariant ring.
The following proposition states that this operator still exists for
invariants of permutation groupoids, albeit missing the important
property of being a $\K[X]^G$-module morphism. In particular, although
$\K[X]^G$ still contains the ring of symmetric polynomials $\sym(X)$,
$R$ is not anymore a $\sym(X)$-module morphism.
\begin{proposition}\label{proposition.groupoids.reynolds}
  There exists an %
  \TODO{unique?}%
  idempotent $R$ in the groupoid algebra $\K.G$ which projects $\K[X]$
  on the invariant ring $\K[X]^G$:
  \begin{displaymath}
    R:=\sum_{A\subset X} \frac{1}{|g\in G, \dom g = A|} \sum_{g\in G,
    \dom g = A} \gr g\ .
  \end{displaymath}
  Furthermore, the four following conditions are equivalent: $R$ is a
  $\sym(X)$-module morphism, $R$
  is a $\K[X]^G$-module morphism, 
  $\ker R$ is a $\sym$-module, and $G$ comes from a permutation group.
\end{proposition}

  


\subsubsection{The chain product}

\FIXME{Double check that the invariant ring is indeed a
  Stanley-Reisner ring, and describe it!}

We now define another product $\star$ on the invariant ring $\K[X]^G$,
called the \emph{chain product}, which preserves a finer grading.
In fact, $(\K[X]^G, \star)$ is a simple realization of the Stanley-Reisner
ring of a suitable poset.  Such rings have been studied
intensively, in particular by Garsia and
Stanton~\cite{Garsia_Stanton.1984} to construct $\sym(X)$-module
generators for the invariant rings of certain permutation groups, and
prove the degree bound for permutation groups
$\beta(G)\leq\binom{|X|}2$ (recall that the \emph{degree bound}
$\beta(A)$ of a finitely generated graded algebra $A$ is the smallest
integer such that $A$ is generated by its elements of degree at most
$\beta(A)$).  This tool is characteristic free: all statements below
actually hold over any ground ring.


Given a subset $S$ of $X$, set $x_S:=\prod_{i\in S} x_i$. By
square-free decomposition, any monomial $x^d$ can be identified
uniquely with a \emph{multichain} $S_1\subset \dots \subset S_k$
of nested subsets of $X$, so that:
\begin{displaymath}
  x^d = x_{S_1} \dots x_{S_k}\ .
\end{displaymath}
We call each $S_k$ a \emph{layer} of $x$. The \emph{fine degree} of the
monomial $x^d$ is the integer vector $(r_1,\dots,r_n)$ where each
$r_i$ counts the (possibly null) number of repetitions of the layer of
size $i$ in $x^d$. The fine degree defines a filtration on $\K[X]$.
%
The \emph{chain product} $\star$ of two monomials $x^d = x_{S_1} \dots
x_{S_k}$ and $x^{d'} = x_{S_1'} \dots x_{S_k'}$ is defined by:
\begin{displaymath}
  x^d \star x^{d'} :=
  \begin{cases}
    x^d x^{d'} &
    \text{ if $\{S_1, \dots, S_k, S'_1,\dots,S'_k\}$ is again a
      multichain of subsets,}\\
    0 & \text{otherwise.}
  \end{cases}
\end{displaymath}
For example, $x_1 \star x_1 = x_1^2$, $x_1 \star x_2=0$, $x_1 x_3^2 \star
x_1x_2x_3^2= x_1^2x_2x_3^4$, and $x_1 x_3^2 \star x_1x_2=0$.

The chain product endows $(\K[X],\star)$ with a second algebra structure
(in fact $(\K[X],\star)$ is isomorphic to the quotient $\K[x_S, S\subset
X] / \{x_Sx_{S'}=0, S\not\subset S'$ and $S'\not\subset S$\}).  It is
also finely graded (fine degrees being added term-by-term).  In fact,
$(\K[X],\star)$ is exactly the associated graded algebra of $\K[X]$ w.r.t.
the fine degree filtration\TODO{Reference?}. Beware that $(\K[X],\star)$
is not an integral domain.

The elementary symmetric functions
\begin{displaymath}
  e_d := \sum_{S\subset X, |S|=d} x_S
\end{displaymath}
are still algebraically independent and generate $(\sym(X)_n,\star)$.
Note that this does not hold for, say, the symmetric powersums. The
following simple fact turns out to be an essential key:
\begin{remark}
  Consider the chain product of a monomial $x_{S_1} \dots x_{S_k}$ by the
  elementary symmetric function $e_d$. It is the sum of all monomials
  $x_{S_1} \dots x_S \dots x_{S_k}$, where $S$ is of size $k$, and
  fits in the chain $S_1\subset\dots \subset S\subset\dots\subset
  S_k$. In particular, if $x_{S_1} \dots x_{S_k}$ readily contains a
  layer $S$ of size $k$, then $x_{S_1} \dots x_{S_k} \star e_k$ is the unique
  monomial obtained by replicating this layer.
\end{remark}

More generally, $(\K[X]^G,\star)$ is a subring of $(\K[X],\star)$. In
particular, $(\K[X]^G,\star)$ is a $\sym(X)$-module. Furthermore, we may
transfer the following algebraic properties from $(\K[X],\star)$ to
$\K[X]^G$, as in the case of permutation
groups~\cite{Garsia_Stanton.1984}.
\begin{proposition}
  \begin{itemize}
  \item[(a)] A family $F$ of finely homogeneous invariants of positive
    degree which generates $(\K[X]^G,\star)$, also generates $\K[X]^G$;
  \item[(b)] $\beta(\K[X]^G,\star) \geq \beta(\K[X]^G)$;
  \item[(c)] A family $F$ of finely homogeneous invariants which
    generates $(\K[X]^G,\star)$ as a $\sym(X)$-module also generates
    $\K[X]^G$ as a $\sym(X)$-module;
  \item[(d)] If $(\K[X]^G,\star)$ is a free $\sym(X)$-module, then so is
    $\K[X]^G$.
  \end{itemize}
\end{proposition}
\begin{proof}
  This is a standard fact about filtrations and associated graded
  connected algebras. The key of the proof is that, if $p$ and $q$ are
  finely homogeneous, the maximal finely homogeneous component of $pq$
  is exactly $p\star q$. (a) and (c) follow by induction over the fine
  grading. Then, (b) follows straightaway from (a), and (d) from (c)
  by a simple Hilbert series argument.
\end{proof}

The converse of (a) and (b) do not hold. In fact, with most
permutation groups, the degree bound $\beta(\K[X]^G,\star)$ is much larger
than $\beta(\K[X]^G)$. We conjecture that the converse of (c) and (d)
hold.  However (d) does not hold anymore in a slightly larger setting
which includes the $r$-quasi-symmetric polynomials of
F.~Hivert~\cite{Hivert.RQSym.2004}, a counter example being
$\qsym^2(X_3)$ (there is an obstruction in the fine Hilbert series).

\begin{theorem}\label{theorem.finiteGeneration}
  The invariant ring $\K[X]^G$ is a finitely generated algebra and
  $\sym(X)$-module, in degree at most $\frac{|X|(|X|+1)}2$. This
  degree bound is tight.
\end{theorem}
Note that, as usual, when $G$ does not act transitively on the
variables, the degree bound can be greatly improved by considering the
elementary symmetric polynomials on each transitive component instead.
\begin{proof}
  The set of orbit sums $o(x_{S_1} \dots x_{S_k})$, where
  $S_1\subsetneq \dots \subsetneq S_k$ is a chain, generate
  $(\K[X]^G,\star)$ as a $(\sym,\star)$-module. This transfers back to
  $\K[X]^G$ and $\sym$.
  
  Note that we may need to consider chains with $S_k=X$; hence the
  degree bound of $\frac{|X|(|X|+1)}{2}$ instead of $\binom{|X|}2$ for
  permutation groups. For an example where the bound is achieved,
  consider the group $G$ made of the identity together with all the
  local bijections of $X=\{1,\dots,n\}$ whose domain is of size at
  most $|X|-1$; then, $\K[X]^G$ is freely generated as a $\sym$-module
  by $1$ and the "staircase" monomials $x_1^{d_1},\dots,x_n^{d_n}$
  with $1 \leq d_i \leq i$.
\end{proof}

\TODO{|
\subsection{(Fine)-Hilbert series}
Polya enumeration?
}

\subsubsection{The Cohen-Macaulay property}

Invariant rings of permutation groups are always Cohen-Macaulay, and
in fact free $\sym(X)$-modules. This follows easily from the
fact that the Reynolds operator is a $\sym(X)$-module morphism.  A
recent and more involved result is that, for all $n$, $\qsym(X_n)$ is
also a free $\sym(X_n)$-module~\cite{Garsia_Wallach.2003}.

As the following example will show, this property does not hold for
all permutation groupoids $G$. Still, $\K[X]^G$ and $(\K[X]^G,\star)$
being finitely generated over $\sym(X)$, they are Cohen-Macaulay if
and only if they are free $\sym(X)$-modules.

\begin{example}
  \label{example.nonCM.groupoid}
  Let $G$ be the permutation groupoid on $\{1,2,3\}$ of
  example~\ref{examples.restrictions}~(a), generated by the local
  bijection $1\mapsto 2$. Then, $G$ is the restriction of a finite
  permutation group whose invariant ring is Cohen-Macaulay. However,
  $\K[X]^G$ itself is not Cohen-Macaulay.  Computing the Hilbert
  series shows right away that this module is not free:
  \begin{displaymath}
    \hilbert_{\K[X]^G} =
    \frac{1}{{\left(1 - Z\right)}^3} - \frac{Z}{1 - Z} =
    \frac{1+Z+2Z^2+2Z^3+Z^4-Z^6}{(1-Z)(1-Z^2)(1-Z^3)}\ .
  \end{displaymath}

  To be more explicit, the transitive components of $G$ being
  $\{1,2\}$ and $\{3\}$, we may replace $\sym(X)$ by
  $R=\sym(x_1,x_2)\otimes \sym(x_3)$, and view $\K[X]^G$ as a finitely
  generated $R$-module. Then, as suggests the Hilbert series,
  \begin{displaymath}
    \hilbert_{\K[X]^G} =
    \frac{1+Z^2+Z^3- Z^4}{(1-Z)^2(1-Z^2)}\ .
  \end{displaymath}
  $\K[X]^G$ is minimally generated as an $R$-module by $(1, x_1 x_3,
  x_1^2 x_2)$, subject to the single relation $x_3.(x_1^2 x_2) = (x_1
  x_2) . (x_1 x_3)$.
  
  Finally, there is no way of choosing the numerator of the Hilbert
  series with non-negative coefficients. Indeed, $\hilbert_{\K[X]^G} =
  \frac{1-Z+2Z^2-Z^3}{(1-Z)^3}$, and the coefficient of highest degree
  in the product of the numerator by
  $\frac{(1-Z^{n_1})(1-Z^{n_2})(1-Z^{n_3})}{(1-Z)^3}$ is always $-1$.
\end{example}

\begin{problem}\label{problem.CM}
  Characterize the permutation groupoids $G$ whose invariant rings
  $\K[X]^G$ (or $(\K[X]^G,\star)$) are Cohen-Macaulay.
\end{problem}

The following theorem is a straightforward extension of a theorem
of~\cite{Garsia_Stanton.1984}.
\begin{theorem}
  $(\K[X]^G,\star)$ is a free $\sym(X)$-module if and only if the
  incidence matrix between generators and maximal chains is
  invertible. In particular, for a set $F$ of finely homogeneous
  invariants whose fine degrees are given by the Hilbert series of
  $\K[X]^G$, the three following conditions are equivalent: $F$ spans
  $\K[X]^G$ as a $\sym(X)$-module, $F$ is a free $\sym(X)$-family, and
  $F$ is a $\sym(X)$-basis of $\K[X]^G$.
\end{theorem}
This readily gives us a necessary condition on the number of
generators.
\begin{corollary}
  If $(K[X]^G,\star)$ is a free $\sym(X)$-module, then it is of rank
  $\frac{|X|!}{|G(X,X)|}$, where $G(X,X)$ is the underlying
  permutation group of $G$.
\end{corollary}

\TODO{In fact, we could do this in a much larger setting; all that is
  required is that two isomorphic monomials have the same shape +
  conditions to have a sym-module, ...  In such settings, the
  invariant ring is always a finite sym-module, but the degree bound
  may be as large as desired.}

\TODO{
  SAGBI-Gröbner w.r.t. the chain product,
  E-R decomposition;
  QSym is not shellable (maybe);
  notion of Break, any relation is obtainable by successive breaks
  from the maximal relation.}

\subsubsection{SAGBI bases}

SAGBI bases (Subalgebra Analog of a Gröbner Bases for Ideals) were
introduced in~\cite{Kapur_Madlener.1989,Robbiano_Sweedler.1990} to
develop an elimination theory in subalgebras of polynomial rings.
Unlike Gröbner bases, not all subalgebras have a finite SAGBI basis,
and it remains a long open problem to characterize those subalgebras
which have a one.
The following theorem states that, as in the case of permutation
groups, invariant rings of permutation groupoids seldom have finite
SAGBI bases. The proof follows the short proof given by the second
author in~\cite{Thiery_Thomasse.SAGBI.2002} for permutation groups,
with some adaptations.  For example $\qsym(X_n)$, represented as a
subring of $\K[X]$, has no finite SAGBI basis whenever $n>1$.  In
particular, $\qsym(X_2)$ becomes the smallest example of finitely
generated algebra which has no finite SAGBI basis (the standard
example being the invariant ring of the alternating group $A_3$).
Still, SAGBI bases and SAGBI-Gröbner bases provide a useful device in
the computational study of invariant rings of permutation
groups~\cite{Thiery.CMGS.2001}, and most likely play the same role
with permutation groupoids.

\begin{theorem}\label{theorem.SAGBI}
  Let $G$ be a permutation groupoid, and $<$ be any admissible term
  order on $\K[X]$. Then, the invariant ring $\K[X]^G$ has a finite
  SAGBI basis w.r.t. $<$ if, and only if, $G$ comes from a permutation
  group generated by reflections (that is transpositions).
\end{theorem}
The following proof is a close variant on the short proof given by the
second author in~\cite{Thiery_Thomasse.SAGBI.2002} in the special case
of permutation groups. For the sake of readability and completeness,
we include it in full here. The key fact is that a submonoid $M$ of
$\N^n$ is finitely generated if, and only if, the convex cone
$C:=\R_+M$ it spans in $\R_+^n$ is finitely generated (that is $C$ is
a \emph{polyhedral cone}). For details, see for
example~\cite[Corollary
2.8]{Bruns_Gubeladze.PolytopesRingsKTheory.2005}. In particular $C$
must be the intersection of finitely many half spaces, and thus closed
for the euclidean topology.
\begin{proof}
  The if-part is easy, a finite SAGBI basis being given by the
  elementary symmetric polynomials in the variables in each
  $G$-transitive components.
  
  Without loss of generality, we may assume $X=\{1,\dots,n\}$ with
  $x_1>\dots>x_n$. Let $M$ be the monoid of initial monomials in
  $\K[X]^G$, seen as a submonoid of $\N^n$, and $C:=\R_+M$ be the
  convex cone it spans in $\R_+^n$.
  
  At this stage, we cannot give an explicit description of $C$, but we
  can construct a convex cone $C'$ which approximates it closely
  enough for our purposes.  By the standard characterization of
  admissible term orders on $\K[X]$, there exists a family of $n$
  linear forms $l=(l_1,\dots,l_n)$ such that $x^d > x^{d'}$ if and
  only if $l(d) >_\lex l(d')$, where we denote by $l(d)$ the $n$-uple
  $l_1(d_1,\dots,d_n),\dots, l_n(d_1,\dots,d_n)$. Given two vectors
  $v$ and $v'$ in $\R_+^n$, we write $v>v'$ if $l(v)>_\lex l(v')$.
  The partial action of $G$ on monomials extends naturally to a
  partial action on $\R_+^n$: whenever the support of
  $v=(v_1,\dots,v_n)$ in contained in the domain of a local bijection
  $f\in G$, $f.v$ is the vector obtained by permuting the non zero
  entries of $v$ according to $f$.  Let $C'$ be the subset of all
  vectors $v$ of $\R_+^n$ such that $v > f.v$ for all $f.v$ in the
  $G$-orbit of $v$. In fact, $C'$ is a convex cone with non empty
  interior (it contains the $n$ linearly independent vectors
  $(1,0,\dots,0), (1,1,0,\dots,0), \dots, (1,\dots,1)$). By
  construction, $M$ consists of the points of $C'$ with integer
  coordinates. It follows that $C\subset C' \overline C$, where
  $\overline C$ is the topological closure of $C$.
  
  Assume now that $M$ is finitely generated. Then, $C$ is a closed
  convex cone, and $C$ and $C'$ simply coincide.
  
  Assume further that $G$ is not generated by transpositions. Then,
  there exists $a<b$ such that the transposition $(a,b)$ is not in
  $G$, while $a$ is in the $G$-orbit of $b$. Choose such a pair $a<b$
  with $b$ minimal. We claim that there is no transposition $(a',b)$
  in $G$ with $a'<b$.  Otherwise, $a$ and $a'$ are in the same
  $G$-orbit, and by minimality of $b$, $(a,a')\in G$; thus,
  $(a,b)=(a,a')(a',b)(a,a')\in G$. Pick $g\in G$ such that $g.b=a$,
  and for $t\geq 0$, define the vector in $\R_+^n$:
  \begin{displaymath}
    u_t :=
    \left(nt,\ (n-1)t,\ \dots,\ (n-b+2)t,\ n-b+1,\ (n-b)t,
      \ \dots,\ t,\ 1\right)\ .
  \end{displaymath}
  Note that $u_1=(n,\dots,1)$ is in $C$, whereas
  $u_0=(0,\dots,0,n-b+1,0,\dots0)$ is not in $C$ because $g.u_0>u_0$.
  
  Take $t$ such that $0<t\leq1$. Then, the vector $u_t$ has no zero
  coefficients, and in particular its $G$-orbit coincides with its
  orbit w.r.t. the underlying permutation group $G(X,X)$.
  Furthermore, the entries of $u_t$ are all distinct, except when
  $t=\frac{n-b}{n-a'}$ for some $a'<b$, in which case the $a'$-th and
  $b$-th entries are equal. Since $(a',b)\not\in G$, the orbit of
  $u_t$ is of size $|G(X,X)|$, and there exists a unique permutation
  $f_t\in G(X,X)$ such that $f_t.u_t$ is in $C$.
  
  Let $t_0=\inf\{t\geq 0, u_t \in C\}$. If $u_{t_0}\not\in C$, then
  $u_{t_0}$ is in the closure of $C$, but not in $C$, a contradiction.
  Otherwise, $u_{t_0}\in C$, and $t_0>0$ because $u_0\not\in C$.  For
  any permutation $f$, $\{f.u_t, t\geq 0\}$ is a half-line; so, $C$
  being convex and closed, $I_f:=\{t, f.u_t\in C\}$ is a closed
  interval $[x_f, y_f]$. For example, $I_\id=[t_0,1]\subsetneq [0,1]$.
  Since the interval $[0,1]$ is the union of all the $I_f$, there
  exists $f\ne\id$ such that $t_0\in I_f$. This contradicts the
  uniqueness of $f_{t_0}$.
\end{proof}

\TODO{Examples: Graph-quasi-symmetric polynomials
}

\subsection{Stability by derivation}

We denote by $\partial_i$ the derivative w.r.t. the variable $x_i$,
and consider the derivation $D := \sum_{i\in X} \partial_i$ on $\K[X]$.
\begin{proposition}\label{proposition.groupoid.derivation}
  Let $G$ be a permutation groupoid. Then $\K[X]^G$ is stable by the
  derivation $D$ if and only if $G$ comes from a permutation group.
  On the other hand, $\K[X]^G$ is always stable w.r.t. the action of
  the \emph{rational Steenrod operators} $S_k := \sum_i x_i^{k+1}
  \partial_i$ for $k\geq 0$ (see~\cite{Hivert_Thiery.SA.2002} for
  details on the rational Steenrod operators).
\end{proposition}
\begin{proof}
  The if part is trivial, since $D$ commutes with the action of the
  symmetric group $\sg_X$ on $\K[X]$. Similarly, the rational Steenrod
  operators always stabilize $\K[X]^G$ because they commute with the
  action of any local bijection on $\K[X]^G$.
  
  Assume now that $\K[X]^G$ is stable by derivation. Let $f:A\mapsto
  B$ be a local bijection such that $A\subsetneq X$, and take $i$ in
  $X\backslash A$. We just need to prove that $f$ extends to a local
  bijection $g$ in $G$ with domain $A\cup \{i\}$.  Applying induction,
  any local bijection in $G$ will then extend to a permutation, as
  desired.
  
  Take a monomial $m$ whose support is $A$ and whose exponents are all
  distinct and at least $2$, and consider the derivation
  $p=D(o(mx_i))$ of the orbitsum of the monomial $mx_i$ in $\K[X]^G$.
  The monomial $m$ occurs in $p$; hence, by invariance of $p$, $f(m)$
  also occurs in $p$, as the derivative of some monomial $g(m x_i)$ in
  the orbit of $mx_i$. By the choice of the exponents of $m$, $f$ and
  $g$ must coincide on $A$, while at the same time $i$ belongs to the
  domain of $g$.
\end{proof}
\begin{example}
  $\qsym(X_2)$ has no graded derivation of degree $-1$.
\end{example}
%
%

\bibliographystyle{alpha}
\bibliography{main}

\def\cprime{$'$}
  \ifx\undefined\allcaps\def\allcaps#1{#1}\fi\ifx\undefined\allcaps\def\allcap%
s#1{#1}\fi
\begin{thebibliography}{NTT04}

\bibitem[BG05]{Bruns_Gubeladze.PolytopesRingsKTheory.2005}
Winfried Bruns and Joseph Gubeladze.
\newblock Polytopes, rings, and k-theory.
\newblock February 2005.

\bibitem[BM04]{Bogvad_Meyer.2003}
Rikard Bogvad and Thomas Meyer.
\newblock On algorithmically checking whether a hilbert series comes from a
  complete intersection.
\newblock Research Reports in Mathematics~5, Department of Mathematics,
  Stockholm University, may 2004.

\bibitem[Cam90]{Cameron.1990}
Peter~J. Cameron.
\newblock {\em Oligomorphic permutation groups}, volume 152 of {\em London
  Mathematical Society Lecture Note Series}.
\newblock Cambridge University Press, Cambridge, 1990.

\bibitem[Cam97]{Cameron.1997}
Peter~J. Cameron.
\newblock The algebra of an age.
\newblock In {\em Model theory of groups and automorphism groups (Blaubeuren,
  1995)}, volume 244 of {\em London Math. Soc. Lecture Note Ser.}, pages
  126--133. Cambridge Univ. Press, Cambridge, 1997.

\bibitem[FP71]{Fraisse_Pouzet.1971}
Roland Fra{\"{\i}}ss{\'e} and Maurice Pouzet.
\newblock Interpr\'etabilit\'e d'une relation pour une cha\^\i ne.
\newblock {\em C. R. Acad. Sci. Paris S\'er. A-B}, 272:A1624--A1627, 1971.

\bibitem[Fra71]{Fraisse.CLM1.1971}
Roland Fra{\"{\i}}ss{\'e}.
\newblock {\em Cours de logique math\'ematique. {T}ome 1: {R}elation et formule
  logique}.
\newblock Gauthier-Villars \'Editeur, Paris, 1971.
\newblock Deuxi\`eme \'edition revue et modifi\'ee, Collection de ``Logique
  Math\'ematique'', S\'erie A, No. 23.

\bibitem[Fra00]{Fraisse.TR.2000}
Roland Fra{\"{\i}}ss{\'e}.
\newblock {\em Theory of relations}, volume 145 of {\em Studies in Logic and
  the Foundations of Mathematics}.
\newblock North-Holland Publishing Co., Amsterdam, revised edition, 2000.
\newblock With an appendix by Norbert Sauer.

\bibitem[GS84]{Garsia_Stanton.1984}
A.~M. Garsia and D.~Stanton.
\newblock Group actions of {S}tanley - {R}eisner rings and invariants of
  permutation groups.
\newblock {\em Adv. in Math.}, 51(2):107--201, 1984.

\bibitem[GW03]{Garsia_Wallach.2003}
A.~M. Garsia and N.~Wallach.
\newblock Qsym over {S}ym is free.
\newblock {\em J. Combin. Theory Ser. A}, 104(2):217--263, 2003.

\bibitem[Hiv04]{Hivert.RQSym.2004}
Florent Hivert.
\newblock Local actions of the symmetric group and generalisations of
  quasi-symmetric functions.
\newblock Preprint, October 2004.

\bibitem[HM88]{Hodkinson_Macpherson.1988}
I.~M. Hodkinson and H.~D. Macpherson.
\newblock Relational structures determined by their finite induced
  substructures.
\newblock {\em J. Symbolic Logic}, 53(1):222--230, 1988.

\bibitem[HT04]{Hivert_Thiery.SA.2002}
F.~Hivert and N.~M. Thi{\'e}ry.
\newblock Deformation of symmetric functions and the rational {S}teenrod
  algebra.
\newblock In {\em Invariant theory in all characteristics}, volume~35 of {\em
  CRM Proc. Lecture Notes}, pages 91--125. Amer. Math. Soc., Providence, RI,
  2004.

\bibitem[KM89]{Kapur_Madlener.1989}
Deepak Kapur and Klaus Madlener.
\newblock A completion procedure for computing a canonical basis for a
  $k$-subalgebra.
\newblock In {\em Computers and mathematics (Cambridge, MA, 1989)}, pages
  1--11. Springer, New York, 1989.

\bibitem[NTT04]{Novelli_Thibon_Thiery.2004}
Jean-Christophe Novelli, Jean-Yves Thibon, and Nicolas~M. Thi{\'e}ry.
\newblock Alg\`ebres de {H}opf de graphes.
\newblock {\em C. R. Math. Acad. Sci. Paris}, 339(9):607--610, 2004.

\bibitem[Pou76]{Pouzet.1976}
Maurice Pouzet.
\newblock Application d'une propri\'et\'e combinatoire des parties d'un
  ensemble aux groupes et aux relations.
\newblock {\em Math. Z.}, 150(2):117--134, 1976.

\bibitem[Pou78]{Pouzet.TR.1978}
Maurice Pouzet.
\newblock {\em Sur la th\'eorie des relations}.
\newblock PhD thesis, Thèse d'état, Universit\'e Claude-Bernard, Lyon 1, 1978.

\bibitem[Pou81]{Pouzet.RPE.1981}
Maurice Pouzet.
\newblock Application de la notion de relation presque-encha\^\i nable au
  d\'enombrement des restrictions finies d'une relation.
\newblock {\em Z. Math. Logik Grundlag. Math.}, 27(4):289--332, 1981.

\bibitem[Pou02]{Pouzet.PR.2002}
Maurice Pouzet.
\newblock The profile of relations.
\newblock Preprint, lecture at Simon Bolivar University, August 2002.

\bibitem[PR86]{Pouzet_Rosenberg.1986}
M.~Pouzet and I.~G. Rosenberg.
\newblock Sperner properties for groups and relations.
\newblock {\em European J. Combin.}, 7(4):349--370, 1986.

\bibitem[RS90]{Robbiano_Sweedler.1990}
Lorenzo Robbiano and Moss Sweedler.
\newblock Subalgebra bases.
\newblock In {\em Commutative algebra (Salvador, 1988)}, pages 61--87.
  Springer, Berlin, 1990.

\bibitem[Thi01]{Thiery.CMGS.2001}
Nicolas~M. Thi{\'e}ry.
\newblock Computing minimal generating sets of invariant rings of permutation
  groups with {SAGBI}-{G}r\"obner basis.
\newblock In {\em Discrete models: combinatorics, computation, and geometry
  (Paris, 2001)}, Discrete Math. Theor. Comput. Sci. Proc., AA, pages 315--328
  (electronic). Maison Inform. Math. Discr\`et. (MIMD), Paris, 2001.

\bibitem[TT04]{Thiery_Thomasse.SAGBI.2002}
N.~M. Thi{\'e}ry and S.~Thomass{\'e}.
\newblock Convex cones and {SAGBI} bases of permutation invariants.
\newblock In {\em Invariant theory in all characteristics}, volume~35 of {\em
  CRM Proc. Lecture Notes}, pages 259--263. Amer. Math. Soc., Providence, RI,
  2004.

\end{thebibliography}

\end{document}